\documentclass{elsart}

\usepackage{amssymb}

\usepackage[english,francais]{babel}

\usepackage[all]{xy}

\newcommand{\ko}{\: , \;}

\newtheorem{theorem}{Theorem}[section]

\newtheorem{e-proposition}[theorem]{Proposition}

\newtheorem{e-definition}[theorem]{Definition\rm}

\newtheorem{theoreme}{Th{\'e}or{\`e}me}[section]
\newtheorem{lemme}[theoreme]{Lemme}

\newtheorem{remarque}{\it Remarque}

\def\og{\leavevmode\raise.3ex\hbox{$\scriptscriptstyle\langle\!\langle$~}}
\def\fg{\leavevmode\raise.3ex\hbox{~$\!\scriptscriptstyle\,\rangle\!\rangle$}}

\begin{document}

\begin{frontmatter}



%
%
\selectlanguage{francais}
\title{Une structure de cat{\'e}gorie de mod{\`e}les de Quillen sur la cat{\'e}gorie des dg-cat{\'e}gories}

\vspace{-1cm}
\selectlanguage{english}
\title{A Quillen model structure on the category of dg categories}



\author{Goncalo Tabuada}
\ead{tabuada@math.jussieu.fr}

\thanks{Soutenu par FCT-Portugal, bourse {\tt SFRH/BD/14035/2003}.}
\address{Universit{\'e} Paris 7 -- Denis Diderot, UMR 7586 du CNRS,
  Case 7012, 2 place Jussieu, 75251 Paris Cedex 05, France}

\begin{abstract}
We construct a cofibrantly generated Quillen model structure on the
category of small differential graded categories.

\vskip 0.5\baselineskip

\selectlanguage{francais}
\noindent{\bf R{\'e}sum{\'e}}
\vskip 0.5\baselineskip
\noindent
Nous construisons une structure de cat{\'e}gorie de mod{\`e}les de Quillen {\`a}
engendrement cofibrant sur la cat{\'e}gorie des petites cat{\'e}gories diff{\'e}rentielles
gradu{\'e}es.

\end{abstract}

\end{frontmatter}

\selectlanguage{english}
\section*{Abridged English version}

Let $k$ be a commutative ring with unit. By a {\em dg category},
we mean a differential graded $k$-category \cite{Keller}
\cite{Drinfeld}. Let $\mbox{{\bf DCAT}}$ be the category of small dg categories. We
remark that it does not have an initial object. Therefore, we consider the
category $\mbox{{\bf DCATp}}$ of small dg categories that have a
distinguished zero object and where the dg functors respect that object.
We have a faithful functor $\mathbf{I}$ from $\mbox{{\bf DCAT}}$ to
$\mbox{{\bf DCATp}}$ which, to a dg category $\mathcal{C}$, associates
the dg category $\mathcal{C}p$ that is obtained from $\mathcal{C}$
by adding a zero object $p$.

We will introduce a cofibrantly generated Quillen model structure on
$\mbox{{\bf DCATp}}$ such that the weak equivalences will be the
quasi-equivalences \cite{Keller}. We will use the recognition theorem
stated in \protect{\cite[2.1.19]{Hovey}}. We now introduce the notations
needed to define the sets $I$ (resp. J) of generating cofibrations
(resp. acyclic generating cofibrations).

Following \cite[3.7.1]{Drinfeld}, we define
$\mathcal{K}$ to be the dg category that has two objects $1$, $2$
and whose morphisms are generated by  $f \in \mathrm{Hom}_{\mathcal{K}}^0 (1,2)$,
$g \in \mathrm{Hom}_{\mathcal{K}}^0 (2,1)$,
$r_1 \in\mathrm{Hom}_{\mathcal{K}}^{-1} (1,1)$,
$r_2 \in \mathrm{Hom}_{\mathcal{K}}^{-1} (2,2)$ and $r_{12}
\in \mathrm{Hom}_{\mathcal{K}}^{-2} (1,2)$ subject to the relations $df=dg=0$,
$dr_1=gf-1_1$, $dr_2 =fg-1_2$ and $dr_{12}=fr_1 - r_2f$. Let
$\mathcal{A}$ be the dg category with one object $3$, such that
$\mathrm{Hom}_{\mathcal{A}}(3,3)=k$. Let $F$ be the dg functor from
$\mathcal{A}$ to $\mathcal{K}$ that sends $3$ to $1$. Let
$\mathcal{B}$ be the dg category with two objects $4$ and $5$ such
that $\mathrm{Hom}_{\mathcal{B}}(4,4)=k$,
$\mathrm{Hom}_{\mathcal{B}}(5,5)=k$,
$\mathrm{Hom}_{\mathcal{B}}(4,5)=0$ and
$\mathrm{Hom}_{\mathcal{B}}(5,4)=0$. Let $n \in \mathbb{Z}$.
Let $S^{n-1}$ be the complex $k[n-1]$ and let $D^n$ be the
mapping cone on the identity of $S^{n-1}$.
We denote by $\mathcal{P}(n)$ the dg category with two objects $6$ and $7$
such that $\mathrm{Hom}_{\mathcal{P}(n)}(6,6)=k$,
$\mathrm{Hom}_{\mathcal{P}(n)}(7,7)=k$,
$\mathrm{Hom}_{\mathcal{P}(n)}(7,6)=0$,
$\mathrm{Hom}_{\mathcal{P}(n)}(6,7)=D^n$ and with composition given by
multiplication. Let $R(n)$ be the dg functor from $\mathcal{B}$ to
$\mathcal{P}(n)$ that sends $4$ to $6$ and $5$ to $7$. Let
$\mathcal{C}(n)$ be the dg category with two objects $8$ and $9$
such that $\mathrm{Hom}_{\mathcal{C}(n)}(8,8)=k$,
$\mathrm{Hom}_{\mathcal{C}(n)}(9,9)=k$,
$\mathrm{Hom}_{\mathcal{C}(n)}(9,8)=0$,
$\mathrm{Hom}_{\mathcal{C}(n)}(8,9)=S^{n-1}$ and composition given by
multiplication. Let $S(n)$ be the dg functor from $\mathcal{C}(n)$
to $\mathcal{P}(n)$ that sends $8$ to $6$, $9$ to $7$ and $S^{n-1}$ to
$D^n$ by the identity on $k$ in degree $n-1$. Finally, let $Q$ be the
dg functor, now in $\mbox{{\bf DCATp}}$, from $\mathcal{O}$, which is the
initial object in $\mbox{{\bf DCATp}}$, to $\mathbf{I}\mathcal{A}$.

\begin{theorem}
If we consider for $\mathcal{C}$ the category $\mbox{{\bf DCATp}}$,
for $W$ the subcategory of quasi-equivalences, for $J$ the functors
$\mathbf{I}F$ and $\mathbf{I}R(n), n\in \mathbb{Z}$, and for $I$ the functors
$Q$ and $\mathbf{I}S(n), n \in \mathbb{Z}$, then the  conditions of the
recognition theorem \protect{\cite[2.1.19]{Hovey}} are fulfilled.
Thus, the category $\mbox{{\bf DCATp}}$ admits a Quillen model structure
whose weak equivalences are the quasi-equivalences.
\end{theorem}

An analogous result for simplicial categories has been obtained in \cite{Bergner}.
Our construction is inspired by \cite{Rezk} and by the construction of
$DG$-quotients in \cite{Drinfeld}. One can easily show that for
this structure, every object is fibrant.

\selectlanguage{francais}
\section{Pr{\'e}liminaires}
Dans toute la suite, $k$ d{\'e}signe un anneau commutatif avec $1$.
Le produit tensoriel $\otimes$ d{\'e}signe toujours le produit tensoriel
sur $k$. Par une {\em dg-cat{\'e}gorie}, nous entendons une $k$-cat{\'e}gorie
diff{\'e}rentielle gradu{\'e}e, voir \cite{Keller} \cite{Drinfeld}.
Soit $\mbox{{\bf DCAT}}$ la cat{\'e}gorie des petites dg-cat{\'e}gories
(elle ne poss{\`e}de pas d'objet initial) et soit $\mbox{{\bf DCATp}}$
la cat{\'e}gorie des petites dg-cat{\'e}gories qui ont un
objet nul sp{\'e}cifi{\'e} et o{\`u} les dg-foncteurs pr{\'e}servent cet objet
nul. On dispose d'un foncteur fid{\`e}le $\mathbf{I}$ de la cat{\'e}gorie $\mbox{{\bf DCAT}}$
vers $\mbox{{\bf DCATp}}$ qui, {\`a} une dg-cat{\'e}gorie $\mathcal{C}$, associe la
dg-cat{\'e}gorie $\mathcal{C}p$ qui s'obtient {\`a} partir de $\mathcal{C}$ en rajoutant un
objet nul $p$. Pour les cat{\'e}gories de mod{\`e}les de Quillen, nous renvoyons {\`a}
\cite{Hovey}. On introduira une structure de cat{\'e}gorie de mod{\`e}les de
Quillen {\`a} engendrement cofibrant dans $\mbox{{\bf DCATp}}$ dont les {\'e}quivalences faibles sont les
quasi-{\'e}quivalences \cite{Keller}. Pour cela, on se servira du th{\'e}or{\`e}me~2.1.19
de \cite{Hovey}.

\section{Th{\'e}or{\`e}me principal}

Suivant \cite[3.7.1]{Drinfeld}, nous d{\'e}finissons $\mathcal{K}$ comme
la dg-cat{\'e}gorie avec deux objets $1$, $2$ et dont
les morphismes sont engendr{\'e}s par $f \in \mathrm{Hom}_{\mathcal{K}}^0 (1,2)$,
$g \in \mathrm{Hom}_{\mathcal{K}}^0 (2,1)$,
$r_1 \in\mathrm{Hom}_{\mathcal{K}}^{-1} (1,1)$,
$r_2 \in \mathrm{Hom}_{\mathcal{K}}^{-1} (2,2)$ et $r_{12}
\in \mathrm{Hom}_{\mathcal{K}}^{-2} (1,2)$ soumis aux relations $df=dg=0$,
$dr_1=gf-1_1$, $dr_2 =fg-1_1$ et $dr_{12}=fr_1 - r_2f$.
$$\xymatrix{
    1 \ar@(ul,dl)[]_{r_1} \ar@/^/[r]^f \ar@/^0.8cm/[r]^{r_{12}} &
    2 \ar@(ur,dr)[]^{r_2} \ar@/^/[l]^g }
$$
Soit $\mathcal{A}$ la dg-cat{\'e}gorie avec un seul object
$3$ et telle que $\mathrm{Hom}_{\mathcal{A}}(3,3)=k$. Soit $F$ le
dg-foncteur de $\mathcal{A}$ vers $\mathcal{K}$ qui envoie $3$
sur $1$.
Soit $\mathcal{B}$ la dg-cat{\'e}gorie avec deux
objects $4$ et $5$ telle que
$
\mathrm{Hom}_{\mathcal{B}}(4,4)=k \ko
\mathrm{Hom}_{\mathcal{B}}(5,5)=k \ko
\mathrm{Hom}_{\mathcal{B}}(4,5)=0  \ko
\mathrm{Hom}_{\mathcal{B}}(5,4)=0
$.
Soit $n \in \mathbb{Z}$. On note $S^{n-1}$ le complexe $k[n-1]$
et $D^n$ le c{\^o}ne sur le morphisme identique de $S^{n-1}$.
On note $\mathcal{P}(n)$ la dg-cat{\'e}gorie avec deux
objets $6$ et $7$ et telle que
$
\mathrm{Hom}_{\mathcal{P}(n)}(6,6)=k \ko
\mathrm{Hom}_{\mathcal{P}(n)}(7,7)=k \ko
\mathrm{Hom}_{\mathcal{P}(n)}(7,6)=0 \ko
\mathrm{Hom}_{\mathcal{P}(n)}(6,7)=D^n
$.
Soit $R(n)$ le dg-foncteur de $\mathcal{B}$ vers
$\mathcal{P}(n)$ qui envoie $4$ sur $6$ et $5$ sur $7$.
On consid{\`e}re la dg-cat{\'e}gorie $\mathcal{C}(n)$ avec deux objects
$8$ et $9$ telle que
$
\mathrm{Hom}_{\mathcal{C}(n)}(8,8)=k \ko
\mathrm{Hom}_{\mathcal{C}(n)}(9,9)=k \ko
\mathrm{Hom}_{\mathcal{C}(n)}(9,8)=0  \ko
\mathrm{Hom}_{\mathcal{C}(n)}(8,9)=S^{n-1}
$.
Soit $S(n)$ le dg-foncteur de $\mathcal{C}(n)$ vers
$\mathcal{P}(n)$ qui envoie $8$ sur $6$, $9$ sur $7$ et $S^{n-1}$ dans $D^n$
par l'identit{\'e} sur $k$ en degr{\'e} $n-1$.
Soit finalement $Q$ le dg-foncteur, maintenant dans $\mbox{{\bf DCATp}}$, de la
cat{\'e}gorie $\mathcal{O}$, qui est l'objet initial dans $\mbox{{\bf DCATp}}$,
vers $\mathbf{I}\mathcal{A}$.

\begin{theoreme} Si on consid{\`e}re pour cat{\'e}gorie $\mathcal{C}$ la cat{\'e}gorie $\mbox{{\bf DCATp}}$,
pour classe $W$ la sous-cat{\'e}gorie de $\mbox{{\bf DCATp}}$ des quasi-{\'e}quivalences,
pour classe $J$ les dg-foncteurs $\mathbf{I}F$ et $\mathbf{I}R(n)$,
$n\in \mathbb{Z}$, et pour classe $I$ les dg-foncteurs $Q$ et $\mathbf{I}S(n)$,
$n\in \mathbb{Z}$, alors les conditions du th{\'e}or{\`e}me \cite[2.1.19]{Hovey} sont satisfaites.
\end{theoreme}

\begin{remarque} Un r{\'e}sultat analogue pour les cat{\'e}gories simpliciales a {\'e}t{\'e}
obtenu dans \cite{Bergner}. Notre construction est inspir{\'e}e par
\cite{Rezk} et par la construction des dg-quotients dans \cite{Drinfeld}.
On peut montrer ais{\'e}ment que pour la structure obtenue, tout objet est fibrant.
\end{remarque}

On observe facilement que les conditions $\mbox{{\it (i)}}$,
$\mbox{{\it (ii)}}$ et $\mbox{{\it (iii)}}$ sont verifi{\'e}es.

\begin{lemme} \label{J-cell-dans-W} On a $J-\mbox{cell}\subseteq W$.
\end{lemme}
Soit $n\in \mathbb{Z}$. Soit $T : \mathbf{I}\mathcal{B} \rightarrow \mathcal{J}$
un morphisme quelconque dans $\mbox{{\bf DCATp}}$. On consid{\`e}re la somme
amalgam{\'e}e suivante
$$
\xymatrix{
\mathbf{I}\mathcal{B} \ar[r]^{T} \ar[d]_{\mathbf{I}R(n)} & \mathcal{J}
\ar[d]^{\mbox{inc}} \\
\mathbf{I}\mathcal{P}(n) \ar[r] & \mathcal{U}
}
$$
dans $\mbox{{\bf DCATp}}$. Il s'agit de v{\'e}rifier que $\mbox{inc}$
est une quasi-{\'e}quivalence.
La cat{\'e}gorie $\mathcal{U}$ s'obtient {\`a} partir de la cat{\'e}gorie
$\mathcal{J}$ en rajoutant un nouveau morphisme
$j$ de $T(4)$ vers $T(5)$ de degr{\'e} $n-1$ et un
nouveau morphisme $l$ de $T(4)$ vers $T(5)$ de
degr{\'e} $n$ tels que $dl=j$. Pour des objets $X$ et $Y$ de $\mathcal{J}$,
on a donc une d{\'e}composition de
$\mathrm{Hom}_{\mathcal{U}}(X,Y)$ en  somme directe de complexes
$$
\mathrm{Hom}_{\mathcal{U}}(X,Y)=\bigoplus_{m \geq 0} \mathrm{Hom}_{\mathcal{U}}^{(m)}(X,Y)
$$
avec
$$
\mathrm{Hom}_{\mathcal{U}}^{(m)}(X,Y)= \underbrace{(T(5),Y) \otimes D^n \otimes
(T(5),T(4)) \otimes D^n \otimes \cdots \otimes D^n
\otimes (X,T(4))}_{m\textrm{ {\scriptsize facteurs} } D^n} \ko
$$
o{\`u} l'on {\'e}crit $(,)$ pour $\mathrm{Hom}_{\mathcal{J}}(,)$.
Puisque le complexe $D^n$ est contractile, l'inclusion
$$
\mathrm{Hom}_{\mathcal{J}}(X,Y) \hookrightarrow \mathrm{Hom}_{\mathcal{U}}(X,Y)
$$
est un quasi-isomorphisme.
Comme le dg-foncteur d'inclusion est l'identit{\'e} au niveau des objets,
c'est une quasi-{\'e}quivalence.
Soit maintenant $N : \mathbf{I}\mathcal{A} \rightarrow \mathcal{L}$ un morphisme
quelconque dans $\mbox{{\bf DCATp}}$.
On consid{\`e}re la somme amalgam{\'e}e suivante
$$
\xymatrix{
\mathbf{I}\mathcal{A} \ar[r]^{N} \ar[d]_{\mathbf{I}F} & \mathcal{L}
\ar[d]^{\mbox{inc}} \\
\mathbf{I}\mathcal{K} \ar[r] & \mathcal{M}
}
$$
dans $\mbox{{\bf DCATp}}$. Il s'agit de montrer que $\mbox{inc}$ est
une quasi-{\'e}quivalence.
La cat{\'e}gorie  $\mathcal{M}$ s'obtient {\`a} partir de la cat{\'e}gorie
$\mathcal{L}$ en rajoutant la cat{\'e}gorie $\mathcal{K}$ {\`a}
$\mathcal{L}$ en identifiant les objets $N(3)$ et $\mathbf{I}F(3)$.
Soit $\mathcal{L}_{0}$ la cat{\'e}gorie $\mathcal{L}$ {\`a} laquelle on rajoute
un morphisme $s$ de $N(3)$ vers un nouvel objet $H$.
Notons $\mathrm{Mod}\mathcal{L}_0$ la cat{\'e}gorie des dg-modules ({\`a} droite) sur
$\mathcal{L}_0$. On consid{\`e}re le plongement de Yoneda
$$
\mathcal{L}_{0} {\hookrightarrow} \mathrm{Mod} \mathcal{L}_0 \ko X \mapsto \widehat{X}.
$$
Soit $\mathcal{L}_{1}$ la sous-cat{\'e}gorie pleine de
$\mathrm{Mod}{\mathcal{L}_{0}}$ dont les objets sont  le c{\^o}ne $C$ sur $\widehat{s}$
et les foncteurs repr{\'e}sentables. Soit
$\mathcal{L}_{2}$ la cat{\'e}gorie obtenue en rajoutant dans $\mathcal{L}_{1}$ un
morphisme $h$ de degr{\'e} $1$ {\`a} l'anneau d'endomorphismes de
$C$  tel que $dh$ est {\'e}gal {\`a} l'identit{\'e} de $C$.
Notre cat{\'e}gorie $\mathcal{M}$ s'identifie naturellement {\`a} la sous-cat{\'e}gorie pleine
de $\mathcal{L}_{2}$  dont les objets sont les images dans $\mathcal{L}_2$
des objets de $\mathcal{L}_0$. Soient $X$ et $Y$ des objets de $\mathcal{L}$.
On a alors une d{\'e}composition de $k$-modules gradu{\'e}es
$$
\mathrm{Hom}_{\mathcal{M}}(X,Y)\simeq
\mathrm{Hom}_{\mathcal{L}_2}(\widehat{X},\widehat{Y}) = \bigoplus_{n\geq 0}
\mathrm{Hom}_{\mathcal{L}_2}^{(n)}(\widehat{X},\widehat{Y}) \ko
$$
o{\`u}
$$
\mathrm{Hom}_{\mathcal{L}_2}^{(n)}(\widehat{X},\widehat{Y})=
\underbrace{\mathrm{Hom}_{\mathcal{L}_1}(C,
\widehat{Y})\otimes S^{2}\otimes \mathrm{Hom}_{\mathcal{L}_1}(C,C)
\otimes S^{2} \otimes \cdots \otimes S^{2}
\otimes \mathrm{Hom}_{\mathcal{L}_1}(\widehat{X},C)}_{n \textrm{ {\scriptsize facteurs} } S^{2}}.
$$
Mais dans cette situation, on n'a pas une somme directe de complexes. Soit
$g_{n+1}\cdotp h \cdotp g_n \cdot h \cdots h \cdotp g_1 \in \mathrm{Hom}_{\mathcal{L}_2}^{(n)}(\widehat{X},\widehat{Y})$.
Comme on a $dh=1$, l'image par $d$ de cet {\'e}l{\'e}ment est {\'e}gale {\`a}
$$
d(g_{n+1})\cdotp h \cdotp g_n \cdotp h  \cdots
 h \cdotp g_1 + \underbrace{(-1)^{\mid g_{n+1} \mid}\cdotp
   g_{n+1}\cdotp 1 \cdotp  g_n \cdotp h  \cdots  h \cdotp g_1}_{(n-1) \textrm{ {\scriptsize facteurs} } h} + \cdots \;\; .
$$
On remarque que, pour tout $m\geq 0$, la somme
$$
\bigoplus_{n \geq 0}^m \mathrm{Hom}_{\mathcal{L}_2}^{(n)}(\widehat{X},\widehat{Y})
$$
est un sous-complexe de $\mathrm{Hom}_{\mathcal{L}_2}(\widehat{X},\widehat{Y})$ et on dispose donc d'une
filtration exhaustive de $\mathrm{Hom}_{\mathcal{L}_2}(\widehat{X},\widehat{Y})$. Le $n$-i{\`e}me
sous-quotient s'identifie {\`a}
$\mathrm{Hom}_{\mathcal{L}_2}^{(n)}(\widehat{X},\widehat{Y})$ et comme le complexe
$\mathrm{Hom}_{\mathcal{L}_1}(\widehat{X},C)$
est contractile l'inclusion
$$
\mathrm{Hom}_{\mathcal{L}}(X,Y) {\hookrightarrow}
\mathrm{Hom}_{\mathcal{M}}(X,Y) \simeq \mathrm{Hom}_{\mathcal{L}_2}(\widehat{X},\widehat{Y})
$$
est un quasi-isomorphisme. Comme $s$ devient un
isomorphisme dans $\mathrm{H}^0(\mathcal{M})$ et que le foncteur d'inclusion est
l'identit{\'e} au niveau des objets, il est bien une quasi-{\'e}quivalence.

D{\'e}montrons maintenant que   $J-\mbox{inj}\cap W = I-\mbox{inj}$.
Pour cela, on consid{\`e}re la classe $\mbox{{\bf Surj}}$ form{\'e}e des foncteurs
$G : \mathcal{H} \rightarrow \mathcal{I}$ dans $\mbox{{\bf DCATp}}$ qui v{\'e}rifient:
$G$ induit une surjection de l'ensemble des objets de $\mathcal{H}$ sur l'ensemble des
objets de $\mathcal{I}$ et $G$ induit des quasi-isomorphismes surjectifs
dans les complexes de morphismes.

\begin{lemme} \label{I-inj-Surj}
On a $I-\mbox{inj} = \mbox{{\bf Surj}}$.
\end{lemme}

Soit $\mathcal{C}$ une cat{\'e}gorie quelconque et $\mathcal{V}$ une classe quelconque
de morphismes dans $\mathcal{C}$. On note $\mathcal{V}-\mbox{drt}$ la classe de morphismes qui
ont la propri{\'e}t{\`e} de rel{\`e}vement {\`a} droite par rapport {\`a} $\mathcal{V}$.
 La classe $Q-\mbox{drt}$ est form{\'e}e des foncteurs qui
  sont surjectifs au niveau des objets. La classe
    $\mathbf{I}S(n)-\mbox{drt}$ est form{\'e}e des foncteurs qui sont des
  quasi-isomorphismes surjectifs au niveau des complexes de
  morphismes. En effet, un carr{\'e} commutatif dans $\mbox{{\bf DCATp}}$
$$
\xymatrix{
\mathbf{I}\mathcal{C}(n) \ar[r]^{D} \ar[d]_{\mathbf{I}S(n)} & \mathcal{H} \ar[d]^{G} \\
\mathbf{I}\mathcal{P}(n) \ar[r]^{E} & \mathcal{I}
}
$$
correspond {\`a} la donn{\'e}e d'un carr{\'e} commutatif dans la cat{\'e}gorie des complexes
$$
\xymatrix{
S^{n-1} \ar@{^{(}->}[d]_{i_n} \ar[r]^-{D} &
  \mathrm{Hom}_{\mathcal{H}}(D(8),D(9)) \ar@{|->}[d]^{G}\\
D^{n} \ar[r]^-{E} & \mathrm{Hom}_{\mathcal{I}}(E(6),E(7))
}
$$
o{\`u} $D(8)$ et  $D(9)$ sont des objets
quelconques dans $\mathcal{H}$. La propri{\'e}t{\'e} r{\'e}sulte de la
caract{\'e}risation des quasi-isomorphismes surjectifs dans la cat{\'e}gorie
des complexes sur $k$. Voir \protect{\cite[2.3.5]{Hovey}}.

\begin{lemme}
\label{J-inj-Surj}
On a $J-\mbox{inj} \cap W = \mbox{{\bf Surj}}$.
\end{lemme}

Montrons l'inclusion $\supseteq$.
Soit $H$ un foncteur de $\mathcal{N}$ vers
$\mathcal{E}$ dans la classe $\mbox{{\bf Surj}}$. Comme $H$
est surjectif au niveau des objets et un quasi-isomorphisme au niveau
des complexes de morphismes, on a $H\in W$.

La classe $R(n)-\mbox{drt}$ est form{\'e}e des foncteurs surjectifs aux
niveau des complexes de morphismes. Il suffit donc de montrer que
$H \in \mathbf{I}F-\mbox{drt}$. La donn{\'e}e d'un carr{\'e} commutatif
$$
\xymatrix{
\mathbf{I}\mathcal{A} \ar[r]^{P} \ar[d]_{\mathbf{I}F} &
  \mathcal{N} \ar[d]^{H} \\
\mathbf{I}\mathcal{K} \ar[r]^{U} & \mathcal{E}
}
$$
correspond {\`a} la donn{\'e}e de la partie inf{\'e}rieure gauche du diagramme
$$
\xymatrix{
P(3) \ar@{|->}[d]_-{H} \ar@{-->}[r]^-{\overline{U(f)}} & D
\ar@{|->}[d]^-{H} \\
U(1) \ar[r]^-{U(f)} & U(2)
}
$$
et {\`a} la donn{\'e} d'une contraction $h$ du c{\^o}ne $C_1$ de $\widehat{U(f)}$ dans $\widehat{\mathcal{E}}$.
Comme $H$ est surjectif au niveau des objets, il
existe $D \in \mathcal{N}$ telle que $H(D)=U(2)$. Le
foncteur $H$ est un quasi-isomorphisme surjectif au niveau
des complexes de morphismes. Donc on peut relever $U(f)$ en
$\overline{U(f)}$.
Dans les cat{\'e}gories des dg-modules respectives, on obtient le
diagramme suivant
$$
\xymatrix{
\widehat{P(3)}\ar@{|->}[d]_-{\widehat{H}} \ar@{-->}[r]^-{\widehat{\overline{U(f)}}} &
\widehat{D} \ar@{|->}[d]^-{\widehat{H}}  \ar[r] &
C_2 \ar@(dr,ur)[]_{h^{\star}} \ar@{|->}[d]^-{\widehat{H}} \\
\widehat{U(1)} \ar[r]^-{\widehat{U(f)}} &
\widehat{U(2)} \ar[r] &  C_1 \ar@(dr,ur)[]_h
}
$$
o{\`u} $C_1$ et $C_2$ d{\'e}signent les c{\^o}nes sur les morphismes
respectifs et  $h$ est la contraction de $C_1$. Comme
$H$ et donc $\widehat{H}$ induisent des
quasi-isomorphismes surjectifs dans les alg{\`e}bres d'endomorphismes,
on peut relever $h$ en une contraction $h^*$ de $C_2$ par application du lemme
\cite[2.3.5]{Hovey} au couple $(h,1)$.

Montrons maintenant l'inclusion $\subseteq$.
Soit $L$ un foncteur de $\mathcal{D}$ vers
$\mathcal{S}$ qui appartient {\`a} $J-\mbox{inj} \cap W$. La classe
$R(n)-\mbox{drt}$ est form{\'e}e des foncteurs surjectifs au niveau
des complexes de morphismes. Comme $L \in W$, il
suffit de montrer que $L$ est surjectif au niveau des
objets.
Soit $E \in \mathcal{S}$ un objet quelconque. Comme $L \in W$,
il existe $C \in \mathcal{D}$ et un morphisme $q \in
\mathrm{Hom}_{\mathcal{S}}(L(C),E)$ qui devient un isomorphisme dans
$\mathrm{H}^{0}(\mathcal{S})$
$$
\xymatrix{
C \ar@{|->}[d]_{L} & \\
L(C) \ar[r]^{q} & E .
}
$$
Ainsi, $q$ est l'image de $f$ par un foncteur de
$\mathbf{I}(\mathcal{K})$ vers $\mathcal{S}$.
Comme on a $L \in J-\mbox{inj}$, on peut relever le
morphisme  $q$ et par cons{\'e}quence l'objet $E$. Le foncteur $L$ est donc bien
surjectif au niveau des  objets.

Nous avons v{\'e}rifi{\'e} que $J-\mbox{cell}\subseteq W$ (lemme~\ref{J-cell-dans-W}) et
que $I-\mbox{inj}$ est {\'e}gal {\`a} $J-\mbox{inj}\cap W$ (lemmes~\ref{I-inj-Surj} et \ref{J-inj-Surj}).
Ces conditions impliquent celles du th{\'e}or{\`e}me de Hovey \cite[2.1.19]{Hovey}.

\label{}


\section*{Remerciements}

Je tiens {\`a} remercier B. Keller pour des conversations utiles.

\end{document}